\documentclass[10pt]{article}

\usepackage{graphicx,amsmath,amssymb,amsthm}

\newtheorem{theorem}{Theorem}[section]
\newtheorem{lemma}[theorem]{Lemma}
\newtheorem{prop}[theorem]{Proposition}
\newtheorem{cor}[theorem]{Corollary}
\theoremstyle{definition}
\newtheorem{defin}[theorem]{Definition}
\newtheorem{rem}[theorem]{Remark}

\newcommand{\Ker}{\mathop{\rm Ker}\nolimits} 
\newcommand{\Ima}{\mathop{\rm Im}\nolimits}   

\newcommand{\f}{\varphi}
\newcommand{\Q}{\mathbb Q} 
\newcommand{\G}{\Gamma} 
 
\newcommand{\C}{{\mathcal C}} 
\newcommand{\ig}{\includegraphics}
\newcommand{\ftimes}{\times_\Phi}

\newcommand{\fdtimes}{\times_{\Phi + d}}

\begin{document}

\centerline{\huge Knight move in chromatic cohomology}
\bigskip

\noindent\parbox[t]{2in}
{\textbf{Michael Chmutov}

Department of Mathematics,

The Ohio State University, 

231 W. 18th Avenue, 

Columbus, Ohio 43210

\texttt{chmutov@mps.ohio-state.edu}
}\qquad
\parbox[t]{2.4in}
{\textbf{Sergei Chmutov}

The Ohio State University, Mansfield,

1680 University Drive,

Mansfield, OH 44906

\texttt{chmutov@math.ohio-state.edu}
}

\bigskip

\noindent
\parbox{2.5in}
{\textbf{Yongwu Rong}

Department of Mathematics,

The George Washington University,

Washington, DC 20052

\texttt{rong@gwu.edu}
}

\newpage

Running head: Knight move in chromatic cohomology

Corresponding author:

Michael Chmutov

8853 Orinda Rd.

Powell, OH 43065

\newpage

\begin{abstract}
In this paper we prove the knight move theorem for the chromatic graph 
cohomologies with rational coefficients introduced by 
L.~Helme-Guizon and Y.~Rong.
Namely, for a connected graph $\G$ with $n$ vertices the only non-trivial 
cohomology groups
$H^{i,n-i}(\G)$, $H^{i,n-i-1}(\G)$ come in isomorphic pairs: 
$H^{i,n-i}(\G)\cong H^{i+1,n-i-2}(\G)$
for $i\geqslant 0$ if $\G$ is non-bipartite, and for $i> 0$ if $\G$ is bipartite.
As a corollary, the ranks of the cohomology groups are determined by the 
chromatic polynomial. At the end, we give an explicit formula for the 
Poincar\'e polynomial in terms of the chromatic polynomial and a 
deletion-contraction formula for the Poincar\'e polynomial. 
\end{abstract}

\section*{Introduction} \label{s:intro}

Recently, motivated by the Khovanov cohomology in knot theory \cite{Kho}, Laure Helme-Guizon and 
Yongwu Rong \cite{HGR1} developed a bigraded cohomology theory for graphs. Its main property is that 
the Euler characteristic with respect to one grading and the Poincar\'e polynomial with 
respect to the other grading give the chromatic polynomial of the graph.

There is a long exact sequence relating the cohomology of a graph with the cohomologies of graphs 
obtained from it by contraction and deletion of an edge. It generalizes the classical 
contraction-deletion rule for the chromatic polynomial. This sequence is an important tool in 
proving various properties of the cohomology (see \cite{HGR1}). In particular, 
for a connected graph $\G$ with $n$ vertices, it allows to prove that the cohomologies are 
concentrated on two diagonals $H^{i,n-i}(\G)$
and $H^{i,n-i-1}(\G)$. We prove that $H^{i,n-i}(\G)$ is isomorphic
to $H^{i+1,n-i-2}(\G)$ (``knight move") for all $i$ with the exception of $i = 0$ for a bipartite 
graph $\G$. An analogous theorem for Khovanov cohomology in knot theory was proved in \cite{Lee}. 
Our proof follows the same idea of considering an additional differential $\Phi$ on the chromatic 
cochain complex of the graph $\G$. This differential anticommutes with the original differential $d$. 
The associated spectral sequence collapses at the term $E_2$, 
and so this term is given by the cohomologies of $\G$ with respect to the differential $\Phi + d$. They 
turn out to be trivial with a small exception (see Theorem \ref{0hompd} for the precise statement). 
On the other hand, the term $E_1$ of the 
spectral sequence is represented by the cohomology of our graph. 
The differential in the term $E_1$ induced by the map $\Phi$ gives the desired isomorphism. Its 
existence implies that the chromatic polynomial of 
a graph determines its cohomology groups with coefficients in a field of characteristic $0$. For 
simplicity we work here with rational coefficients. The constuction of the cohomology is based on the 
algebra $\Q[x]/(x^2)$ and reflects its properties. In general, for chromatic cohomologies based 
on other algebras \cite{HGR2}, the cohomology groups are supported on more than two diagonals, and 
we do not expect them to be determined by the chromatic polynomial.

It would be interesting to adapt the recent techniques of spanning trees \cite{Weh, CK} and the Karoubi 
envelopes \cite{BNM} from the Khovanov cohomology theory to the chromatic cohomology. This might lead to 
another proof of our key Theorem \ref{0hompd} and to a deeper understanding of combinatorics.

This work was motivated by computer calculations of the chromatic homology by M.~Chmutov, which 
revealed certain patterns in the Betti numbers. Part of this work has been completed  
during the Summer'05 VIGRE working group ``Knot Theory and Combinatorics" at the Ohio State University 
funded by NSF, grant DMS-0135308. Y.~Rong was partially supported by the NSF grant DMS-0513918. 
The authors would like to thank L.~Helme-Guizon and J.~Przytycki for numerous discussions,
S.~Duzhin and anonymous referees for valuable comments.
\section{Definitions and preliminary results} \label{s:def}

\begin{figure} [bt]
\ig[width=300pt] {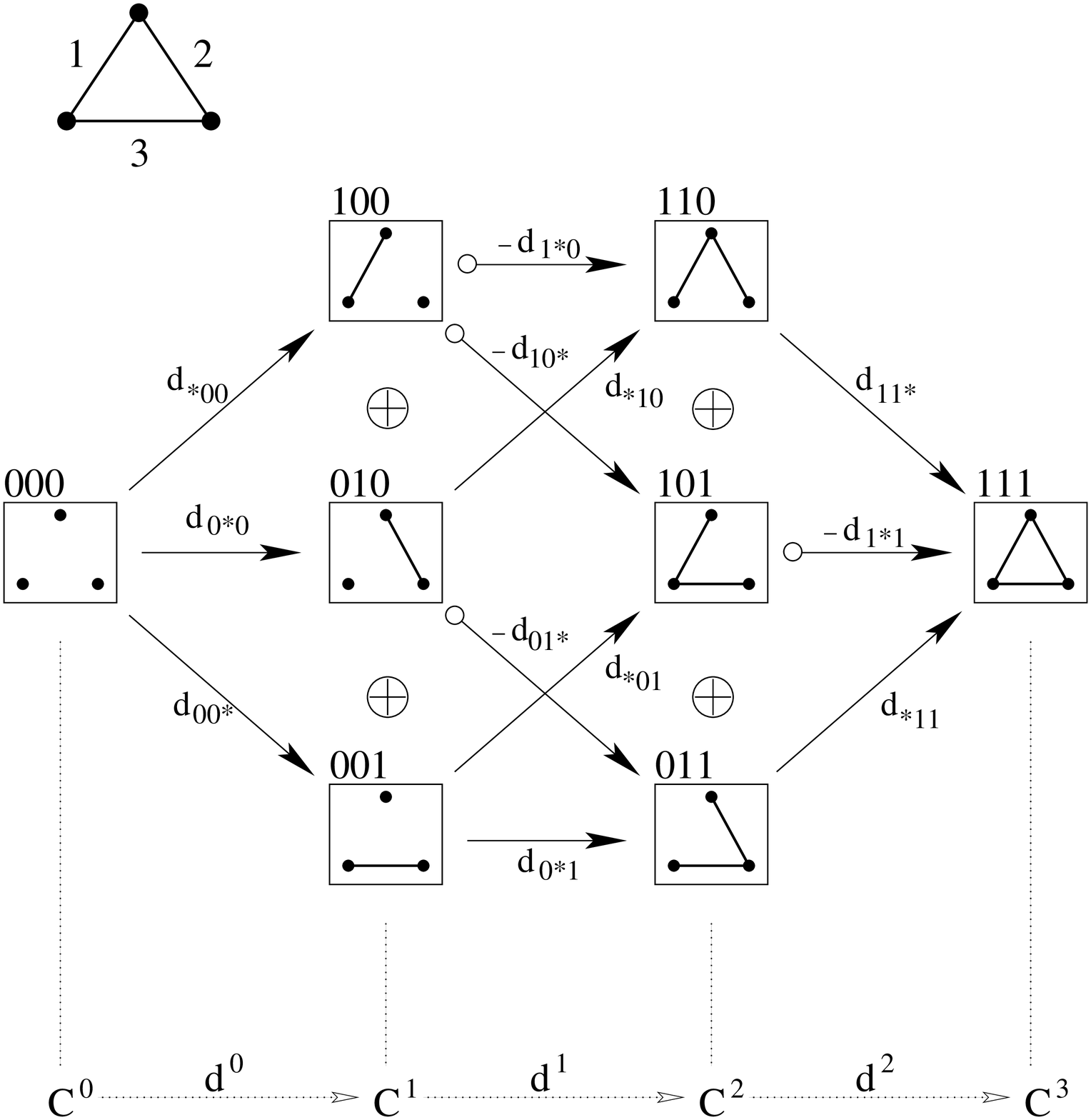}
\caption{Chromatic cochain complex.}
\label{bnp}
\end{figure}

For a graph $\G$ with ordered edges, a {\it state} $s$ is a 
spanning subgraph of $\G$, that is a 
subgraph of $\G$ containing all the vertices of $\G$ and a subset of the 
edges. The number of edges in a state is called its {\it dimension}. An 
{\it enhanced state} $S$ is a state whose connected components are colored in 
two colors: $x$ and $1$. The number of connected components colored in $x$ 
is called the {\it degree} of the enhanced state. The {\it cochain group} 
$C^{i,j}$ is defined to be the real vector space spanned by all enhanced 
states of dimension $i$ and degree $j$. These notions are illustrated on 
Figure \ref{bnp} similar to Bar-Natan's
\cite{BN}.
Here every square box represents a vector space spanned by all enhanced 
states with the indicated underlying state. The direct sum of these vector 
spaces located in the $i$-th column gives the cochain group 
$C^i=\bigoplus\limits_j C^{i,j}$.
The boxes are labeled by strings of $0$'s and $1$'s which encode the 
edges participating in the corresponding states. To turn the cochain groups 
into a {\it cochain complex} we define a differential 
$d^{i,j}:C^{i,j}\to C^{i+1,j}$.
On a vector space corresponding to a given state (box) the differential can 
be defined as adding an edge to the corresponding state in all possible 
ways, and then coloring the connected components of the obtained state 
according to the following rule. Suppose we have an enhanced state $S$ with 
an underlying state $s$ and we are adding an edge $e$. Then, if the number 
of connected components is not 
changed, we preserve the same coloring of connected components of the new 
state $s\cup \{e\}$. If $e$ connects two different connected components of $s$, 
then the color of the new component of $s\cup \{e\}$ is defined by the multiplication
$$1\times 1:=1,\qquad 1\times x:=x,\qquad x\times 1:=x,\qquad x\times x:=0\ .$$
In the last case, the enhanced state is mapped to zero. 
In cases where the number of edges of $s$ whose order is less than that of $e$ is odd, 
we take the target enhanced state $S\cup \{e\}$ with the coefficient $-1$. 
These are shown in the picture above by arrows with little circles at their tails. 
See \cite{HGR1} for more details as well as for a proof of the main property 
$d^{i+1,j}\circ d^{i,j}=0$ converting our cochain groups into a bigraded cochain 
complex $C^{*,*}(\G)$. We call its cohomology groups the {\it chromatic cohomology} 
of the graph $\G$:
$$H^{i,j}(\G)\ :=\ 
\frac{\Ker\left(d: C^{i,j}(\G)\to C^{i+1,j}(\G)\right)}{
      \Ima\left(d: C^{i-1,j}(\G)\to C^{i,j}(\G)\right)}\ .
$$

The vector space corresponding to a single vertex without edges is isomorphic to 
the algebra of truncated polynomials 
${\mathcal A}:=\Q[x]/(x^2)$. Using this algebra we can think about a box space 
of an arbitrary graph $\G$ as a tensor power of the algebra 
${\mathcal A}$ whose tensor factors are in one-to-one correspondence
with the connected components of the state. Then our multiplication rule for 
the differential turns out to be the multiplication operation 
${\mathcal A}\otimes {\mathcal A}\to {\mathcal A}$
in the algebra ${\mathcal A}$. This approach allows to generalize the definition 
of the chromatic cohomology to an arbitrary algebra 
${\mathcal A}$ (see \cite{HGR2,Pr} for further development of this approach).
In particular, the cohomology of $\G$ with respect to the differential 
$\Phi + d$ that we will study later can be understood as the cohomology associated 
with the algebra
$\C = \Q[x]/(x^2-1)$.
 
The following facts are known (see \cite{HGR1, HGPR})
\begin{itemize}
\item If $P_\G(\lambda)$ denotes the chromatic polynomial of $\G$ then
$$\displaystyle P_\G(1 + q) = \sum_{i, j} (-1)^iq^j\dim (H^{i,j}(\G))\ .$$
\item For a graph $\G$, having an edge $e$, let $\G - e$ and $\G / e$ denote the graphs 
obtained from $\G$ by deletion and contraction of $e$, respectively. Then, there exists 
a long exact sequence (for any $j$)
$$
0\rightarrow H^{0,j}(\G)
   \rightarrow H^{0,j}(\G-e)
   \rightarrow H^{0,j}(\G/e)
   \rightarrow H^{1,j}(\G)
   \rightarrow\dots
$$
\item For a graph $\Gamma$ with $n\geqslant 2$ vertices, $H^{i,*}(\G)=0$ if $i>n-2$.
\item For a graph $\Gamma$ with $n$ vertices and $k$ connected components the 
cohomologies are concentrated on $k+1$ diagonals:
 $H^{i,j}(\G)=0$ unless $n-k\leqslant i+j\leqslant n$.
\item For a loopless connected graph $\Gamma$ with $n$ vertices only the 
following $0$-cohomologies are nontrivial:
$H^{0,n}(\G)\cong H^{0,n-1}(\G)\cong \Q$ for bipartite graphs, and $H^{0,n}(\G)\cong \Q$ for 
non-bipartite graphs. 
\end{itemize}
\section{The differentials $\Phi$ and $\Phi + d$}
\begin{defin}
Define a map $\Phi: C^i(\G) \rightarrow C^{i+1}(\G)$ in the same way as  $d$ except using the algebra ${\mathcal B}$:
$$
1  \ftimes  1  =  0,
1  \ftimes  x  =  0,
x  \ftimes  1  =  0,
x  \ftimes  x  =  1.
$$
\end{defin}
Note that the algebra ${\mathcal B}$ is not unital.

\begin{prop}
The map $\Phi$ is a differential, i.e. $\Phi^2 = 0$.
\end{prop}

This is a censequence of the algebra structure on ${\mathcal B}$. See \cite{HGR2, Pr} for details.

\begin{prop} 
The cochain bicomplex $C^{*,*}(\G)$ with differentials $d$ and $\Phi$ is independent of the 
ordering of edges.
\end{prop}

\begin{proof}
The proof is the same as that of \cite[Theorem 14]{HGR1}, where an isomorphism $f$ between 
chain complexes with different edge orderings was constructed. The isomorphism $f$ also 
commutes with $\Phi$.
\end{proof}

\begin{prop} 
$\Phi\circ d + d\circ\Phi = 0$.
\end{prop}

\begin{proof} 

The map $\Phi + d: C^i(\G) \rightarrow C^{i+1}(\G)$ can be described in the same way as $d$
except using the algebra $\C=\Q[x]/(x^2 - 1)$:
$$
1  \fdtimes  1  =  1,\qquad
1  \fdtimes  x  =  x,\qquad
x  \fdtimes  1  =  x,\qquad
x  \fdtimes  x  =  1.
$$
Of course, $\C$ is isomorphic to $\Q^2$ as a vector space.
Therefore $\Phi + d$ is a differential. But then
$$0 = (\Phi + d)^2 = \Phi^2 + (\Phi\circ d + d\circ\Phi) + d^2,$$
and so $\Phi\circ d + d\circ\Phi = 0$.

\end{proof}
\section{Cohomology of $\Phi+d$}

Consider the cohomologies $H^i_{\Phi + d}$ with respect 
to the differential $\Phi+d$. These cohomologies are not graded 
spaces anymore, since $\Phi + d$ does not preserve the grading.
Instead, they have a natural filtration which is 
preserved by $\Phi + d$. We will discuss this in more detail at 
the end of the section.

\begin{rem}
The map $\Phi + d$ commutes with the maps from 
the short exact sequence of complexes associated with $\G / e$, $\G$, and $\G - e$,
in the same way as $d$ does. Therefore is a long exact sequence of the new 
cohomology groups:
$$
  0\rightarrow H^0_{\Phi + d}(\G)
   \rightarrow H^0_{\Phi + d}(\G - e)
   \stackrel{\gamma}{\rightarrow} H^0_{\Phi + d}(\G / e)
   \rightarrow H^1_{\Phi + d}(\G)
   \rightarrow H^1_{\Phi + d}(\G - e)
   \rightarrow \cdots
$$
As in \cite{HGR1}, all cohomology groups of a graph containing a loop are 
trivial and that multiple edges of a graph can be replaced by single ones 
without altering the cohomology groups.
\end{rem}

\begin{rem}\label{basic states}
Let $a_0=\frac{1}{2}(x+1), a_1=\frac{1}{2}(x-1)$. Then 
$\{ a_0,a_1\}$ forms a new basis for the algebra $\C$ with
$$ a_0  \fdtimes  a_0  =  a_0,  \qquad
a_1  \fdtimes  a_1  =  -a_1, \qquad
a_0  \fdtimes  a_1  =  a_1 \fdtimes  a_0  =  0.
$$
Under the new basis, the enhanced states $S$ are now spanning graphs
$s$ whose components are colored with $a_0$ and $a_1$. If
$\Gamma$ is bipartite, we can split its vertex set
$V(\Gamma)=V_0\cup V_1$ into the two parts $V_0$ and $V_1$.
This gives two particular enhanced states in dimension 0: $S_0$
(resp. $S_1$) is the coloring on $V(\Gamma)$ with each vertex in
$V_i$ (resp. $V_{1-i}$) colored by $a_i$. The property that $a_0
\fdtimes  a_1  =  a_1 \fdtimes a_0 = 0$ immediately implies
\end{rem}

\begin{lemma}
$(\Phi + d)(S_0)=(\Phi + d)(S_1)=0$, in other words, $S_0$ and $S_1$ are both cocycles
in $C_{\Phi + d}^0(\G)$.
\end{lemma}

In fact, $\{ S_0, S_1\}$ forms a basis of $H_{\Phi + d}^0(\G)$ as the next
theorem states.

\begin{theorem}
\label{0hompd}
Let $\G$ be a connected graph.\\
{\rm 1.} If\ \ $\G$ is not bipartite, then $H^i_{\Phi + d}(\G)=0$ for all $i$.\\
{\rm 2.} If\ \ $\G$ is bipartite, then $H^i_{\Phi + d}(\G)=0$ for
all $i>0$,  $H^0_{\Phi + d}(\G)\cong \Q^2$ with basis  
$\{ S_0, S_1\}$ described above.
\end{theorem}

\begin{proof}
We induct on $m$, the number of edges in $\G$.

If $m=0$, $\G$ consists of one vertex and no edges. We have
$H^0(\G)\cong \C \cong \Q^2$, and $H^i(\Gamma)=0$ for all
$i>0$.

Suppose that the theorem is true for all connected graphs with less 
than $m$ edges. Let $\Gamma$ be a connected graph with
$m$ edges. We consider two cases.

Case A. $\Gamma$ is a tree. Let $e$ be a pendant edge. By
\cite[Proposition 3.4]{HGR2}, $H^i_{\Phi + d}(\G) \cong H^i_{\Phi +
d}(\G / e)\otimes {\mathcal A}'$ where 
${\mathcal A}'=\langle x\rangle \cong \Q$ is a subspace in $\C$ spanned by $x$. 
The graph
$\G/e$ is a tree with one less edge and therefore we can apply
induction. It follows that $H^0_{\Phi + d}(\G)\cong \C \cong \Q^2$,
and $H^i_{\Phi + d}(\G)=0$ for all $i>0$.

Case B. $\G$ is not a tree. It must contain an edge $e$ that is
not a bridge. Consider the long exact sequence
$$
  0\rightarrow H^0_{\Phi + d}(\G)
   \rightarrow H^0_{\Phi + d}(\G - e)
   \stackrel{\gamma}{\rightarrow} H^0_{\Phi + d}(\G / e)
   \rightarrow H^1_{\Phi + d}(\G)
   \rightarrow H^1_{\Phi + d}(\G - e) = 0
$$
Here the last term $H^1_{\Phi + d}(\G - e) = 0$ since $\Gamma - e$
satisfies the induction hypothesis. By Lemma \ref{triple} below, we
have three subcases.

{\it Subcase 1.} All three graphs $\G, \G - e, \G/e$ are
non bipartite. In this case, we have $H^0_{\Phi + d}(\G - e)= H^0_{\Phi + d}(\G /
e)=0$ by the induction hypothesis. Hence 
$H^0_{\Phi + d}(\G)=H^1_{\Phi + d}(\G)=0$.

{\it Subcase 2.} The graphs $\G$ and $\G -e$ are bipartite, while $\G/ e$
is not. In this case, we have 
$H^0_{\Phi + d}(\G / e)=0$, $H^0_{\Phi + d}(\G - e)\cong \Q ^2$. 
It follows that $H^0_{\Phi + d}(\G )\cong \Q ^2$,
$H^i_{\Phi + d}(\G )=0$ for all $i>0$.

{\it Subcase 3.} The graphs $\G-e$ and $\G /e$ are bipartite, while $\G$ 
is not. In this case, we have 
$H^0_{\Phi + d}(\G - e)\cong H^0_{\Phi + d}(\G / e) \cong \Q
^2$. The spaces have bases $\{ S_0(\G -e), S_1(\G -e)\}$ and
$\{ S_0(\G /e ),S_1(\G /e )\}$, respectively. The two endpoints of $e$ in 
$\G-e$ must be labeled by the same color, for otherwise $\G /e$ would
not be bipartite. The property $a_i  \fdtimes  a_i  =  \pm a_i$ then
implies that the connecting map $\gamma$ sends $S_i(\G -e)$ to 
$\pm S_i(\G/e)$ for $i=0, 1$. 
Indeed, $\gamma$ acts on an enhanced state $S$ by inserting 
the edge $e$ and adjusting the coloring according to the multiplication 
rule.
Therefore $\gamma$ is an isomorphism. It follows
that $H^0_{\Phi + d}(\G )= H^1_{\Phi + d}(\G )=0$.

\end{proof}

\begin{lemma}\label{triple}
Let $\G$ be a connected graph and let $e\in E(\G)$ be an edge that
is not a bridge. Then the possible bipartiteness
of the triple $\G, \G-e, \G/e$ is shown in the following table
($\checkmark$ stands for bipartite graphs, while -- stands for non bipartite ones).

\begin{center}\begin{tabular}{c||c|c|c}
Case & $\G$ & $\G-e$ & $\G/e$ \\
\hline\hline
1 & -- & -- & -- \\
2 & \checkmark & \checkmark & -- \\
3 & -- & \checkmark & \checkmark
\end{tabular}\end{center}
\end{lemma}

\begin{proof}
By our assumption, the three graphs $\G, \G-e, \G/e$ are all
connected. Thus by definition, $\G$ (resp. $\G-e, \G/e$) is
bipartite if and only if $P_{\G}(2)=2$ (resp. $P_{\G -e}(2)=2,
P_{\G/e}(2)=2$). On the other hand, the contraction-deletion rule
says $P_{\G -e}(2)=P_{\G}(2)+P_{\G/e}(2)$.

If $\G -e$ is not bipartite, then $P_{\G -e}(2)=0$ which implies that 
$P_{\G}(2)=P_{\G/e}(2)=0$ and therefore $\G$ and $\G/e$ are both non
bipartite. However, this is exactly the first case.

If $\G-e$ is bipartite, then $P_{\G -e}(2)=2$ which implies that either
$P_{\G}(2)=2$ and $P_{\G/e}(2)=0$, or $P_{\G}(2)=0$ and
$P_{\G/e}(2)=2$. The first possibility yields case $2$ while second
possibility yields case $3$.
\end{proof}

\subsection*{Filtered Cohomology}

The cochain bicomplex $C^{*,*}(\G)$ has a natural filtration
$C^{i,\leqslant j}(\G):=\!\bigoplus\limits_{k=0}^j C^{i,k}(\G)$;
$$C^{i,\leqslant 0}(\G)\subseteq C^{i,\leqslant 1}(\G)\subseteq
C^{i,\leqslant 2}(\G)\subseteq \dots \subseteq C^{i,\leqslant n}(\G)
=C^i(\G)\ .
$$

The differential $\Phi + d$ preserves this filtration, because $d$ has 
bidegree $(1,0)$ and $\Phi$ has bidegree $(1,-2)$. So one can talk 
about the cohomology groups $H_{\Phi+d}(C^{*,\leqslant j}(\G))$.

Since for every $j$ there is an embedding of complexes 
$C^{*,\leqslant j}(\G)\subseteq C^*(\G)$,
we have the corresponding homomorphism of cohomology groups
$$H_{\Phi+d}(C^{*,\leqslant j}(\G))\to H_{\Phi+d}(C^*(\G))\ .$$ 
We denote the image of this homomorphism by 
$H_{\Phi+d}^{*,\leqslant j}(\G)$.
Thus we have a filtration
$$H_{\Phi+d}^{i,\leqslant 0}(\G)\subseteq 
  H_{\Phi+d}^{i,\leqslant 1}(\G)\subseteq 
  H_{\Phi+d}^{i,\leqslant 2}(\G)\subseteq \dots \subseteq 
  H_{\Phi+d}^{i,\leqslant n}(\G) = H_{\Phi+d}^i(\G)\ .
$$

It is easy to see that the basic cocycles $S_0$, $S_1$ constructed in Remark
\ref{basic states} both have degree $n$. So they belong to 
$C^{0, \leqslant n}_{\Phi + d}(\G)$, however their difference $S_0-S_1$ belongs to
$C^{0, \leqslant n - 1}_{\Phi + d}(\G)$.
As a direct consequence of this and Theorem \ref{0hompd}
we have the following description of the filtered cohomologies.

\begin{cor}
\label{filtdescr}
For a connected graph $\G$ with $n$ vertices:\vspace{8pt}\\
1. if $\G$ is not bipartite, then  $H^{i, \leqslant j}_{\Phi + d}(\G) = 0$, for 
      all $i$ and $j$;\vspace{8pt}\\
2. if $\G$ is bipartite, then $H^{i, \leqslant j}_{\Phi + d}(\G) = 0$, for 
      $i\geqslant 1$, and\vspace{8pt}\\
\hspace*{20pt}$H^{0, \leqslant n}_{\Phi + d}(\G) = H^0_{\Phi + d}(\G) \cong \Q^2$,\qquad 
$H^{0, \leqslant n - 1}_{\Phi + d}(\G)\cong \Q$,\qquad 
$H^{0, \leqslant n - 2}_{\Phi + d}(\G) = 0$.
%
\end{cor}

\section{Knight Move}
\bigskip

Since $\Phi$ anticommutes with $d$, it acts (``knight move'') on the
$d$-cohomology groups $\Phi\colon H^{i,j}(\G)\to H^{i+1,j-2}(\G)$.
$$\ig[width=340pt] {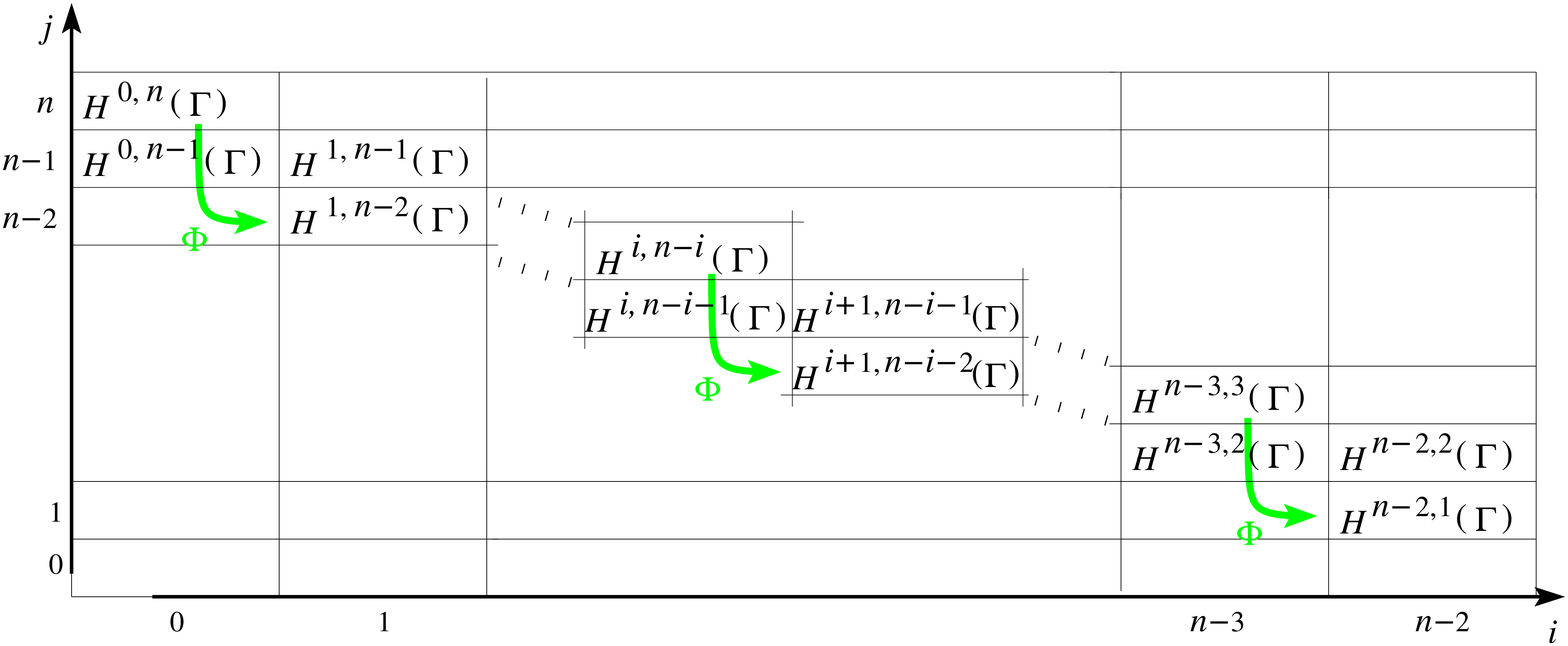}$$

The following theorem is a version of a standard theorem for spectral 
sequences in homological algebra. 

\bigskip

\begin{theorem}[Knight Move]\label{knightmove}\hspace*{-3pt}
Let $\G$ be a connected graph with $n$ vertices. Then there is an isomorphism 
$\f$ of the following quotient spaces
$$\f:\ 
H_{\Phi+d}^{i,\leqslant j}(\G) / H_{\Phi+d}^{i,\leqslant j-1}(\G)
\ \cong\ 
\frac{\Ker\left(\Phi: H^{i,j}(\G)\to H^{i+1,j-2}(\G)\right)}{
      \Ima\left(\Phi: H^{i-1,j+2}(\G)\to H^{i,j}(\G)\right)}.
$$
\end{theorem}

\begin{proof}
For every bicomplex $\{C^{*,*}(\G),d,\Phi\}$ there is a standard 
way to associate a spectral sequence (see \cite[p.47]{MC}, except our 
indices $i,j$ are different). Our original complex
$\{C^{*,*}(\G),d\}$ is the $E_0$ term of the spectral sequence. Our cohomology 
groups $\{H^{*,*}(\G),\Phi\}$ together with the differential
$\Phi: H^{i,j}(\G)\to H^{i+1,j-2}(\G)$ form the term $E_1$. Its 
cohomology groups, which are on the right-hand side of our isomorphism $\f$:
$$E_2^{i,j}(\G) =
 \frac{\Ker\left(\Phi: H^{i,j}(\G)\to H^{i+1,j-2}(\G)\right)}{
      \Ima\left(\Phi: H^{i-1,j+2}(\G)\to H^{i,j}(\G)\right)}
$$
form the term $E_2$.
Its differential has bidegree $(1,-4)$. When the cohomology groups are 
concentrated on two diagonals it is too ``long", so it is zero.
Therefore the spectral sequence collapses at the term $E_2$. 
In other words, $E_{\infty}=E_2$. 
The standard theorem (\cite[Theorem 2.15]{MC}) claims that the 
spectral sequence converges to the bigraded vector space associated 
with the filtration
$$H_{\Phi+d}^{i,\leqslant 0}(\G)\subseteq 
  H_{\Phi+d}^{i,\leqslant 1}(\G)\subseteq 
  H_{\Phi+d}^{i,\leqslant 2}(\G)\subseteq \dots \subseteq 
  H_{\Phi+d}^{i,\leqslant n-i}(\G) = H_{\Phi+d}^i(\G)\ .
$$
That is, the spaces $E_2^{i,j}=E_{\infty}^{i,j}$ are isomorphic to the
corresponding quotient spaces of the filtration 
$H_{\Phi+d}^{i,\leqslant j}(\G) / H_{\Phi+d}^{i,\leqslant j-1}(\G)$.
So our theorem is a direct consequence of this general theorem for 
spectral sequences.
\end{proof}

\begin{cor}\label{knightisom}
For a connected graph $\G$ with $n$ vertices:\\
1. If $\G$ is not bipartite, $\Phi: H^{i,n - i}(\G)\to H^{i+1,n-i-2}(\G)$ is an isomorphism for \\ \mbox{\qquad all $i$.}\\
2. If $\G$ is bipartite, $\Phi: H^{i,n-i}(\G)\to H^{i+1,n-i-2}(\G)$ is an isomorphism for all \\ \mbox{\qquad}
                                                          $i\geqslant 1$;
the map $\Phi: H^{0,n}(\G)\to H^{1,n - 2}(\G)$ has one-dimensional kernel.
\end{cor}

This is a direct consequence of Corollary \ref{filtdescr} and Theorem \ref{knightmove}.

\section{Applications}

In this section we give three applications of the knight move theorem. 
The first one is a computation of the $1$-dimensional homologies 
of a connected graph. The second is an expression for the Poincar\'e 
Polynomial in terms of the chromatic polynomial. This result shows 
that the ranks of the homologies carry no new information as compared to 
the chromatic polynomial. The last theorem shows that the long exact 
sequence of homologies splits into a collection of short exact sequences
starting with the term $H^1(\G/e)$.

The key observation here is that for a connected graph $\G$ with $n$ vertices, 
$$H_{\Phi+d}^{i,\leqslant j}(\G) / H_{\Phi+d}^{i,\leqslant j-1}(\G)$$
is almost always $0$. In fact, it is only non-zero for bipartite $\G$,
for $i = 0$, and for $j=n$ or $n-1$ (Corollary \ref{filtdescr}). In this 
particular case, $H_{\Phi+d}^{0,\leqslant n}(\G) \cong \Q^2$, 
while $H_{\Phi+d}^{0,\leqslant n-1}(\G) \cong \Q$, and 
$H_{\Phi+d}^{0,\leqslant n}(\G) / H_{\Phi+d}^{0,\leqslant n-1}(\G) \cong \Q$.
So, according to the knight move theorem, there is an isomorphism almost 
everywhere between the two diagonals. 

\subsection{Homologies of dimension $1$}

Recall that a {\it simple graph} is a graph that does not have more than one edge 
between any two vertices and no edge starts and ends at the same vertex.

\begin{theorem}\label{t:1hom}
Let $\G$ be a connected and simple graph with $n$ vertices and $m$ edges. 
Then 
$$\begin{array}{lll}
H^{1,n - 1}(\G) = \Q^{m - n + 1} & H^{1, n - 2}(\G) = 0    & \text{, if $\G$ is bipartite}\\
H^{1,n - 1}(\G) = \Q^{m - n} & H^{1, n - 2}(\G)\cong\Q & \text{, otherwise}
\end{array}
$$
\end{theorem}

\begin{proof}
Since the cohomologies are non-zero only on the two diagonals, the only degrees we 
need to worry about are $n-1$ and $n-2$. First, let us find the cohomologies of 
degree $n-1$ (which does not require the knight move theorem).

Notice that $d(C^{1, n-1}(\G)) = 0$ since $C^{2, n-1}(\G) = 0$ for a simple graph. Hence 
$$\begin{array}{lcl}
{\rm dim}(H^{1, n-1}(\G))\! &=\!& \dim(C^{1, n-1}(\G)) -      \dim(d(C^{0, n-1}(\G))) \\
                          &=\!& \dim(C^{1, n-1}(\G)) - 
                             (\dim(C^{0, n-1}(\G)) -      \dim(H^{0, n-1}(\G)))\\
                          &=\!& \dim(C^{1, n-1}(\G)) - (n - \dim(H^{0, n-1}(\G)))\\
                          &=\!& m - n + \dim(H^{0, n-1}(\G)).
\end{array}$$
According to \cite[Theorem 39]{HG}, 
$$
H^{0, j}(\G) \cong\Q\text{ \ for\quad} \begin{cases}
j = n\text{ and }j = n-1  & \text{, if $\G$ is bipartite}\\
j = n                     & \text{, otherwise}
\end{cases}
$$

Hence, if $\G$ is bipartite, ${\rm dim}(H^{1, n-1}(\G)) = m - n + 1$ and if not, then 
${\rm dim}(H^{1, n-1}(\G)) = m - n$.

Now let us calculate the cohomologies of degree $n-2$. If $\G$ is not bipartite, then,
by the knight move theorem, $\Phi: H^{0, n}(\G)\rightarrow H^{1, n-2}(\G)$ is an 
isomorphism. But, again using \cite[Theorem 39]{HG}, $H^{0, n}(\G)\cong \Q$. 
Therefore  $H^{1, n-2}(\G)\cong \Q$.

Next consider the case when $\G$ is bipartite. From the knight move theorem, we know 
that $\Ker(\Phi:H^{0, n}(\G) \rightarrow H^{1, n - 2}(\G)) \cong \Q$. However, 
$H^{0, n}(\G)$ is itself isomorphic to $\Q$. So, 
$\Ima(\Phi:H^{0, n}(\G) \rightarrow H^{1, n - 2}(\G)) = 0.$
Applying the knight move to the next step gives us 
$$0 = \frac{\Ker(\Phi:H^{1, n-2}(\G) \rightarrow H^{2, n - 4}(\G))}
          {\Ima(\Phi:H^{0, n}(\G) \rightarrow H^{1, n - 2}(\G))} \cong H^{1, n-2}(\G),$$
since $H^{2, n - 4}(\G)$ lies off the two diagonals. So $\G$ has no $1$-dimensional 
cohomologies in degree $n-2$.
\end{proof}

\subsection{The Poincar\'e polynomial and the chromatic polynomial}

The Poincar\'e polynomial $R_\G(t,q) := \sum_{i, j} t^iq^j\dim (H^{i,j}(\G))$ of a 
connected graph $\G$ splits into two 
homogeneous parts, $R^n_\G(t,q)$ and $R^{n-1}_\G(t,q)$, since all cohomologies are 
concentrated on two diagonals.

\begin{theorem} For a connected graph $\G$
$$R^n_\G = 
\begin{cases}
\frac{q^2}{t} R^{n-1}_\G + q^n\frac{t - q}{t}  & \text{, if $\G$ is bipartite}\vspace{10pt}\\
\frac{q^2}{t} R^{n-1}_\G                       & \text{, otherwise}
\end{cases}$$
\end{theorem}

\begin{proof}
Suppose $\G$ is not bipartite. Then $H^{0, n - 1} (\G) = 0$. So the conclusion is just a 
consequence of Corollary \ref{knightisom}. 

Now suppose $\G$ is bipartite. Then the situation is similar to the previous case, except
\par 1. The polynomial for the lower diagonal gains a term $q^{n-1}$ coming from $H^{0, n-1}$, and
\par 2. The polynomial for the upper diagonal gains a previously unaccounted term $q^n$ which has been 
        mapped to $0$ by $\Phi$. 
These two complications accout for the two additional terms in the expression for $R^n_\G(t,q)$.
\end{proof}

\begin{rem}
We know that $R_\G(t,q) = R^n_\G(t,q) + R^{n-1}_\G(t,q)$.

If $\G$ is not bipartite, this gives $R_\G(t,q) = \left(1 + \frac{q^2}{t}\right) R^{n - 1}_\G(t,q)$. 
Plugging in $t = -1$ gives $P_\G(1 + q) = (1 - q^2) R^{n - 1}_\G(-1, q)$, 
or $\frac{P_\G(1 + q)}{(1 - q^2)} = R^{n-1}_\G(-1, q)$. However, since 
each term in $R^{n-1}_\G(t,q)$ has degree $n-1$, knowlege of 
$R^{n-1}_\G(-1,q)$ is sufficient to fully determine $R^{n-1}_\G(t,q)$. 
Therefore $R_\G(t,q)$ is determined by $P_\G(1 + q)$.

If $\G$ is bipartite, this gives 
$\displaystyle\frac{P_\G(1 + q) - q^n - q^{n + 1}}{(1 - q^2)} = R^{n-1}_\G(-1, q)$. 
This also means that $R_\G(t, q)$ is determined by $P_\G(1 + q)$.
\end{rem}

\begin{cor} \label{poi-chrom}
For a connected graph $\G$
$$R_\G (t,q) =
\begin{cases}\displaystyle
(-1)^{n-1} t^n\frac{t+q^2}{t^2-q^2}P_\G\left(\frac{t-q}{t}\right) + \frac{q^n(t + 1)}{t + q}
                                & \text{, if $\G$ is bipartite} \vspace{10pt}\\
\displaystyle
(-1)^{n-1} t^n\frac{t+q^2}{t^2-q^2}P_\G\left(\frac{t-q}{t}\right)        
                                & \text{, otherwise}
\end{cases}$$
\end{cor}

\begin{proof}
This is a direct result of the calculation described in the remark.
\end{proof}

\subsection{Deletion-contraction formula for the Poincar\'e polynomial}

Using the bipartiteness table of the triple
$\G$, $\G-e$, $\G/e$ from Lemma \ref{triple} we get the following theorem.

\begin{theorem}
Let $\G$ be a simple, connected graph with $n$ vertices, and let $e$ be an edge that is 
not a bridge. Then in cases $1$ and $2$ (when $\G / e$ is not bipartite)
$$R_\G(t, q) = R_{\G - e}(t, q) + t R_{\G / e}(t, q)\ ,$$
while in case 3 (when $\G$ is not bipartite but $\G - e$ and $\G / e$ are)
$$R_\G(t, q) = R_{\G - e}(t, q) + t R_{\G / e}(t, q) - q^{n - 1}(t + 1)\ .$$
\end{theorem}

The theorem follows from Corollary \ref{poi-chrom} and the deletion-contraction formula 
for the chromatic polynomial $P_\G(\lambda)$.

It implies the following relation between the dimensions of the homology spaces
$$\dim (H^{i,j}(\G)) = \dim (H^{i,j}(\G-e)) + \dim (H^{i-1,j}(\G/e))$$
for all $i$ and $j$ in cases $1$ and $2$, and for $(i,j)\neq (0,n-1)$ or $(1,n-1)$ in case 3.
For the exceptional values of $(i,j)$ in case 3 we have
$$\begin{array}{c||c|c|c}
(i,j) & \dim (H^{i,j}(\G)) & \dim (H^{i,j}(\G-e)) & \dim (H^{i-1,j}(\G/e)) \\
\hline\hline
(0,n-1) & 0 & 1 & 0 \\
(1,n-1) & m-n & m-n & 1
\end{array}
$$
where  $m$ is the number of edges of $\G$.

This relation between the dimensions gives the following splitting of the long exact sequence into short ones.

\begin{prop}
Let $\G$ be a simple, connected graph with $n$ vertices, and let $e$ be an edge that is not a 
bridge. Then the connecting homomorphisms $\varphi$ in the long exact sequence
$$\begin{array}{r}
0\rightarrow H^{0,j}(\G)
   \rightarrow H^{0,j}(\G-e)
   \stackrel{\varphi}{\rightarrow} H^{0,j}(\G/e)
   \rightarrow H^{1,j}(\G)
   \rightarrow H^{1,j}(\G-e) \hspace{40pt}\vspace{10pt}\\
   \stackrel{\varphi}{\rightarrow} H^{1,j}(\G/e)
   \rightarrow H^{2,j}(\G)
   \rightarrow H^{2,j}(\G-e)
   \stackrel{\varphi}{\rightarrow} H^{2,j}(\G/e)
   \rightarrow H^{3,j}(\G)
   \rightarrow\dots
\end{array}
$$
are identically $0$ unless $i=0$, $j = n-1$, $\G$ is not bipartite, 
while $\G - e$ and $\G/e$ are bipartite.

In the exceptional case where $\G$ is not bipartite, 
while $\G - e$ and $\G/e$ are bipartite the connection map
$\varphi\colon H^{0,n-1}(\G-e)\to H^{0,n-1}(\G/e)$
is an isomorphism of one-dimensional vector spaces.

\end{prop}

\end{document}